\documentclass[a4paper]{amsart}

\begin{document}

\title{An improved local wellposedness result for the modified KdV-equation}

\author{Axel Gr\"unrock \\Fachbereich C: Mathematik/Naturwissenschaften\\ Bergische Universit\"at Wuppertal \\ Gau{\ss}stra{\ss}e 20 \\ D-42097 Wuppertal \\ Germany \\ e-mail Axel.Gruenrock@math.uni-wuppertal.de}

\date{}

\maketitle

\newcommand{\R}{\mathbb R}
\newcommand{\h}[2]{\mbox{$ \widehat{H^{#1}_{#2}}$}}
\newcommand{\hh}[3]{\mbox{$ \widehat{H^{#1}_{#2, #3}}$}}
\newcommand{\n}[2]{\mbox{$ \| #1 \| _{ #2} $}}
\newcommand{\x}{\mbox{$X^r_{s,b}$}}
\newcommand{\xx}{\mbox{$X_{s,b}$}}
\newcommand{\X}[3]{\mbox{$X^{#1}_{#2,#3}$}}
\newcommand{\XX}[2]{\mbox{$X_{#1,#2}$}}
\newcommand{\q}[2]{\mbox{$ {\| #1 \|}^2_{#2} $}}

\pagestyle{plain}
\rule{\textwidth}{0.5pt}

\newtheorem{lemma}{Lemma}
\newtheorem{kor}{Corollary}
\newtheorem{satz}{Theorem}
\newtheorem{prop}{Proposition}

\begin{abstract} The Cauchy problem for the modified KdV-equation
\[u_t + u_{xxx} = (u^3)_x, \hspace{2cm} u(0)=u_0\]
is shown to be locally wellposed for data $u_0$ in the space $\h{r}{s}(\R)$ defined by the norm
\[\n{u_0}{\h{r}{s}} := \n{\langle \xi \rangle ^s\widehat{u_0}}{L^{r'}_{\xi}},\]
provided $\frac{4}{3}< r \le 2$, $s \ge \frac{1}{2} - \frac{1}{2r}$. For $r=2$ this coincides with the best possible result on the $H^s$-scale due to Kenig, Ponce and Vega. The proof uses an appropriate variant of the Fourier restriction norm method and linear as well as bilinear estimates for the solutions of the Airy-equation.
\end{abstract}

\section{Introduction} The Cauchy problem for the modified KdV- (mKdV-)equation
\begin{equation}\label{1}
u_t + u_{xxx} = (u^3)_x, \hspace{2cm} u(0)=u_0
\end{equation}
is known to be locally well posed for data $u_0$ in the classical Sobolev spaces $H^s (\R)$ if $s \ge \frac{1}{4}$, and ill posed in the sense that the mapping data upon solution is no longer uniformly continuous, if $s< \frac{1}{4}$. Both, the positive and the negative result, are due to Kenig, Ponce and Vega, see Thm. 2.4 in \cite{KPV93}, respectively Thm. 1.3 in \cite{KPV01}. The standard scaling argument here suggests local wellposedness for $s>- \frac{1}{2}$ and is thus misleading in this case. A very similar situation arises for the semilinear Schr\"odinger equation in one space dimension
\[iu_t + u_{xx} = |u|^{\alpha}u, \hspace{2cm}0< \alpha < 4 ,\]
for which the Cauchy problem is known to be locally (and globally) well posed, if $s \ge 0$ (see \cite{CW90} and the references therein), and ill posed in the sense mentioned above, if $s<0$ (\cite{KPV01}, Thm. 1.1). Again, the scaling argument is misleading since it suggests LWP on the $H^s$-scale for $s > \frac{1}{2} - \frac{2}{\alpha}$. In the Schr\"odinger context it was suggested by Vega and other authors to leave the $H^s$-scale in order to prove local and global wellposedness results for data \emph{not} belonging to $L^2$ anymore, see \cite{VV01}, where the case $\alpha = 2$ is considered, as well as \cite{CVV01}, where the Fourier transform of the data is assumed to be in a weak $L^p$-space for - in the onedimensional case - some $p \in (2,4)$. The crucial linear estimate in the onedimensional part of Thm. 2 of \cite{CVV01} is
\begin{equation}\label{2}
\n{e^{it\partial ^2} u_0}{L^p_{xt}} \le c \n{\widehat{u_0}}{L^{r'}_{\xi}},
\end{equation}
($p=3r$, $\frac{1}{r} + \frac{1}{r'}=1$, $\frac{4}{3} < r \le 2$, which goes back to Fefferman and Stein (\cite{F70}). In \cite{CVV01} the authors restrict themselves to nonlinearities with $\alpha > \frac{8}{3}$, in order to show global and scattering results for small data. This restriction is no longer necessary, if one is interested in \emph{local} wellposedness only. Then, in the cubic (i. e. \\ $\alpha = 2$) case (ignoring the "weak"-refinement for the sake of simplicity), the following result can be easily derived by the aid of the Fefferman-Stein-estimate (\ref{2}):

\begin{prop}\label{p1} Let $\frac{4}{3} < r \le 2$ and $\widehat{u_0} \in L^{r'}(\R)$. Then the Cauchy problem
\begin{equation}\label{3}
iu_t + u_{xx} = |u|^{2}u, \hspace{2cm}u(0)=u_0
\end{equation}
is locally well posed.
\end{prop}
Unfortunately, allthough it is essentially contained in the arguments of \cite{CVV01}, the above proposition is not mentioned explicitely in that paper; on the other hand I cannot see either how to conclude it directly from the local results in \cite{VV01}, which are certainly much deeper. So let me sketch the proof briefly, this will give some light on what follows:

\vspace{0.3cm}

The contraction mapping principle is applied to the integral equation corresponding to the Cauchy problem (\ref{3}) in the closed ball of radius $2c \n{\widehat{u_0}}{L^{r'}}$ in the space
\[C([0,T], \widehat{L^r}) \cap L^{3r}([0,T],L^{3r}_x),\]
where $c$ is the largest constant in the subsequent estimates and
\[\widehat{L^r}:= \{ f \in \mathcal{S}'(\R): \n{f}{\widehat{L^r}}:=\n{\hat{f}}{L^{r'}}<\infty\}.\]
Then the linear part can be estimated by
\[\sup_{t \in [0,T]} \n{e^{-it \xi ^2}\widehat{u_0}}{L^{r'}} = \n{\widehat{u_0}}{L^{r'}},\]
which is trivial, and by
\[\n{e^{it\partial ^2} u_0}{L^{3r}([0,T],L^{3r}_x)} \le c \n{\widehat{u_0}}{L^{r'}},\]
where the Fefferman-Stein-estimate (\ref{2}) comes in. Now using Minkowsky's integral inequality we obtain for the nonlinear part $G(t)= \int_0^t e^{i(t-s)\partial ^2}|u(s)|^2u(s)ds$
\begin{eqnarray*}
&& \sup_{t \in [0,T]}\n{G(t)}{\widehat{L^r}} + \n{G}{L^{3r}([0,T],L^{3r}_x)}\\
& \le & \int_0^T \n{|u(s)|^2u(s)}{\widehat{L^r}}+\n{e^{i(t-s)\partial ^2}|u(s)|^2u(s)}{L^{3r}([0,T],L^{3r}_x)} ds \\
& \le & c \int_0^T \n{u^3(s)}{L^r} ds,
\end{eqnarray*}
where in the last step the Hausdorff-Young inequality and the estimate (\ref{2}) were applied. Finally, H\"older's inequality gives the upper bound
\[... \le c T^{\frac{1}{r'}} \|u\|^3_{L^{3r}([0,T],L^{3r}_x)}.\]
In a similar manner the corresponding difference estimate can be derived. Now choose $T$ sufficiently small and the job is done.

\vspace{0.3cm}

The aim of the present paper is to show a result corresponding to the above proposition for the modified KdV-equation. More precisely: We shall prove the local wellposedness of the Cauchy problem (\ref{1}) for data $u_0 \in \h{r}{s}(\R)$, $\frac{4}{3}<r \le 2$, $s \ge s(r):= \frac{1}{2} - \frac{1}{2r}$, where the space $\h{r}{s}(\R)$ is defined by the norm
\[\n{u_0}{\h{r}{s}} = \n{\widehat{J^s u_0}}{L^{r'}},\]
here $J^s$ is the Bessel potential operator of order $-s$. For $r=2$ this coincides with the LWP-result in \cite{KPV93} mentioned above. In contrast to the Schr\"odinger case, we can \emph{lower} the bound on $s$ with decreasing $r$. It should be mentioned that from the scaling point of view the spaces $\h{r}{s}$ behave like the Bessel potential spaces $H^{s,r}$ (which are embedded in $\h{r}{s}$ for $r \le 2$ by Hausdorff-Young) and like $H^{\sigma}$, if $s-\frac{1}{r} + \frac{1}{2} = \sigma$. For the admissible values of $s$ and $r$ in our result we have $s-\frac{1}{r} + \frac{1}{2} \ge s(r)-\frac{1}{r} + \frac{1}{2} = 1 - \frac{3}{2r}> - \frac{1}{8}$, which certainly can be seen as an improvement (compared with the $r=2$-case), but which is still far away from the bound $s-\frac{1}{r} + \frac{1}{2} > - \frac{1}{2}$ suggested by the scaling argument.

\vspace{0.3cm}

To prove the result we use an appropriate variant of the Fourier restriction norm method introduced by Bourgain in \cite{B93}. This variant is described in a more general setting in section 2. A central argument in our proof is the analogue of the Fefferman-Stein-estimate (\ref{2}) for the Airy-equation, which is shown in section 3 and in which we obtain a gain of almost $\frac{1}{4}$ fractional derivative. This of course is not enough to compensate the "loss" of a whole derivative in the nonlinearity. So the linear estimate has to be supplemented by a bilinear one exhibiting a larger gain of derivatives. This bilinear estimate is also shown - together with some corollaries - in section 3. As the proof shows, it is closely related to the Airy-version of the Fefferman-Stein-estimate mentioned above. Finally, the fourth section is devoted to the proof of the crucial nonlinear estimate.

\section{A variant of Bourgain's method}

For a smooth phase function $\phi: \R ^n \rightarrow \R$ of polynomial growth we define the function spaces
\[\x:= \{f \in \mathcal{S'}(\R ^{n+1}): \n{f}{\x}< \infty\},\]
where $s,b \in \R$, $1 \le r \le \infty$, $\frac{1}{r} + \frac{1}{r'}=1$ and
\[\n{f}{\x}:= \left(\int d \xi d \tau \langle \xi \rangle^{sr'}\langle \tau - \phi(\xi)\rangle^{br'} |\hat{f}(\xi , \tau)|^{r'} \right) ^{\frac{1}{r'}},\]
with the usual modifikation for $r=1$. For $r=2$ these spaces coincide with those introduced by Bourgain in \cite{B93} in the study of initial value problems. So in this case we shall omit the index $r$.

\vspace{0.3cm}

In this section we follow closely the exposition for $r=2$ in \cite{G96}, chapter 3.

\subsection{Elementary properties}

The $\x$-spaces are Banach spaces. For $1<r$ they are separable and contain the Schwartz class $\mathcal{S}(\R ^{n+1})$ as a dense subspace. For $1<r<\infty$ the mapping
\[\Phi: \X{r'}{-s}{-b} \rightarrow (\x)',\,\,\,\,g \mapsto \Phi(g),\]
defined by
\[\Phi(g)[f]:= \int d \xi d \tau \hat{f}(\xi , \tau)\overline{\hat{g}}(\xi , \tau)\]
is isometric, antilinear and onto; so - with respect to the inner product on $L^2_{xt}$ - the dual space $(\x)'$ of $\x$ can be identified with $\X{r'}{-s}{-b}$. From Thm. 5.5.3 in \cite{BL} we can conclude for $s_0, s_1, b_0, b_1 \in \R$, $1<r_0, r_1 \le \infty$, $\theta \in [0,1]$ and
\[s = (1- \theta )s_0 + \theta s_1, \hspace{1cm}b = (1- \theta )b_0 + \theta b_1,\hspace{1cm} \frac{1}{r}=\frac{1- \theta}{r_0} + \frac{\theta}{r_1}\]
that
\[(\X{r_0}{s_0}{b_0},\X{r_1}{s_1}{b_1})_{[\theta]}= \x ,\]
where $[\theta]$ denotes the complex interpolation method. A simple application of H\"older's inequality gives the following continuous embedding:

\begin{equation}\label{4}
\X{r_1}{s_1}{b_1} \subset \X{r_0}{s_0}{b_0} ,
\end{equation}
provided $r_1 \le r_0$, $s_1 - \frac{n}{r_1} > s_0 - \frac{n}{r_0}$ and $b_1 - \frac{1}{r_1} > b_0 - \frac{1}{r_0}$.

\vspace{0.3cm}

The connection between the $\x$-norms and the evolution operators ${U_{\phi}}(t)= e^{it \phi (D)}$, $t \in \R$, is the same as in the $r=2$-case: Defining $\hh{r}{s}{b}$ by the norm
\[\n{f}{\hh{r}{s}{b}}:= \left(\int d \xi d \tau \langle \xi \rangle^{sr'}\langle \tau \rangle^{br'} |\hat{f}(\xi , \tau)|^{r'} \right) ^{\frac{1}{r'}}\]
and using
\[\widehat{{U_{\phi}}(-\cdot) f} (\xi , \tau) = \widehat{f}(\xi, \tau + \phi(\xi))\]
we see that
\begin{eqnarray*}
\n{{U_{\phi}}(-\cdot)f}{\hh{r}{s}{b}} & = & \left(\int d \xi d \tau \langle \xi \rangle^{sr'}\langle \tau \rangle^{br'} |\widehat{{U_{\phi}}(-\cdot)f}(\xi , \tau)|^{r'} \right) ^{\frac{1}{r'}} \\
& = & \left(\int d \xi d \tau \langle \xi \rangle^{sr'}\langle \tau - \phi(\xi)\rangle^{br'} |\hat{f}(\xi , \tau)|^{r'} \right) ^{\frac{1}{r'}} = \n{f}{\x}.
\end{eqnarray*}

Using this, another type of embeddings can be derived from space-time-estimates for the group $(U_{\phi}(t))_{t \in \R}$ :

\begin{lemma}\label{l1} Assume $Y \subset \mathcal{S'}(\R ^{n+1})$ to be a Banach space being stable under multiplication with $L_t^{\infty}$, that is
\[ \n{\psi u}{Y} \leq c \n{\psi }{L_t^{\infty}} \n{\psi u}{Y} \,\,\, \mbox{for all} \,\,\,\psi \in L_t^{\infty}, \,\,\, u \in Y ,\]
such that the inequality
\begin{equation}\label{51}
\n{U_{\phi} u_0}{Y} \leq c \n{u_0 }{\widehat{L^r}}
\end{equation}
holds for all $u_0 \in \widehat{L^r}$. Then for all $ b > \frac{1}{r}$ the estimate
\[\n{ u}{Y} \leq c \n{u}{\X{r}{0}{b}}\]
is valid with a constant c depending only on $b$.
\end{lemma}

Proof: With $g := \mathcal{F} _t U_{\phi}(-\cdot)u$, where $\mathcal{F} _t$ denotes the Fourier transform in the time variable only, we can write
\begin{eqnarray*}
u(t) &=& U_{\phi}(t)U_{\phi}(-t)u(t) \\
&=& c U_{\phi}(t) \int e^{it \tau} (\mathcal{F} _t U_{\phi}(-\cdot)u)(\tau) d\tau \\
&=& c \int e^{it \tau}U_{\phi}(t)g(\tau) d\tau.
\end{eqnarray*}
Now using Minkowski's inequality, the stability assumption on $Y$ and (\ref{51}) we obtain
\begin{eqnarray*}
\n{ u}{Y} &\leq & c \int d \tau \n{U_{\phi}g(\tau)}{Y} \\
&\leq & c \int d \tau \langle \tau \rangle^{-b} (\langle \tau \rangle^{b}\n{g(\tau)}{\widehat{L^r}}).
\end{eqnarray*}
Finally H\"older's inequality is applied and the proof is complete.
$\hfill \Box$

\vspace{0.3cm}

The above lemma can easily be generalized to multilinear estimates. We shall not make use of this here, except in the case $r=2$, where this is well known, see e. g. Proposition 3.5 in \cite{KS01}. A simple but important consequence of Lemma 1 is the embedding
\begin{equation}\label{52}
\x \subset C(\R, \h{r}{s})
\end{equation}
for all $b > \frac{1}{r}$, which follows from the fact that the evolution operators $(U_{\phi}(t))_{t \in \R}$ form a group of isometries on the space $\h{r}{s}$. This will guarantee the persistence property of the solution in our application.

\subsection{Linear estimates}

The identity $\n{{U_{\phi}}(-\cdot)f}{\hh{r}{s}{b}}=\n{f}{\x}$ immediately gives the necessary estimate for the solutions of the homogeneous linear equation:

\begin{equation}\label{5}
\n{\psi {U_{\phi}} u_0}{\x} = \n{\psi u_0}{\hh{r}{s}{b}}= \n{\psi}{\h{r}{b}}\n{u_0}{\h{r}{s}}=: c_{\psi}\n{u_0}{\h{r}{s}}
\end{equation}
for any $C_0^{\infty}$-function $\psi$ of the time variable only and for any $u_0 \in \h{r}{s}$. The next aim is to obtain an estimate for the  solution of the inhomogeneous linear equation
\[\partial _t v -i \phi (D) v=F,\hspace{1cm}v(0)=0,\]
which is given by
\[v(t)= \int_0^t {U_{\phi}} (t-t')F(t')dt'=:{U_{\phi}} _{*R} F(t). \]
For this purpose let $\psi \in C^{\infty}_0 $ with $supp(\psi) \subset (-2,2)$ and, for $0< \delta \leq 1$, $\psi_{\delta}(t)= \psi(\frac{t}{\delta})$.
\begin{lemma}[Estimate for the homogeneous linear equation]\label{l2}
Assume $1<r<\infty$ and $b'+1 \ge b \ge 0 \ge b' > -\frac{1}{r'}$. Then
\[\n{\psi_{\delta}{U_{\phi}} _{*R} F}{\x}\leq c \delta^{1+b'-b}\n{F}{\X{r}{s}{b'}}.\]
\end{lemma}

Proof: First we show for $Kg(t):= \psi_{\delta}(t) \int_0^t g(t')dt'$ that
\begin{equation}\label{6}
\n{Kg}{\h{r}{b}} \leq c \delta^{1+b'-b} \n{g}{\h{r}{b'}} .
\end{equation}
Here the function $g$ is - at first - assumed to depend on the time variable only. Writing
\[ \int_0^t g(t')dt' = c \int_{- \infty}^{\infty} \frac{\exp{(it \tau)}-1}{i \tau}\widehat{ g}(\tau)d \tau \]
we have $Kg(t)= I + II + III$ with
\begin{eqnarray*}
I & = & \psi_{\delta} \sum_{k \geq 1} \frac{t^k}{k !} \int_{|\tau| \delta \leq 1} (i \tau)^{k-1 }\widehat{ g}(\tau)d \tau \\
II & = & -\psi_{\delta} \int_{|\tau| \delta \geq 1} (i \tau)^{-1}\widehat{ g}(\tau)d \tau \\
III & = & \psi_{\delta}\int_{|\tau| \delta \geq 1} (i \tau)^{-1} \exp{(it \tau)}\widehat{ g}(\tau)d \tau .
\end{eqnarray*}
The first contribution can be estimated for $1 \geq b \geq 0 \geq b'$ as follows:
\[\n{I}{\h{r}{b}} \leq \sum_{k \geq 1} \frac{1}{k !} \n{t^k \psi_{\delta}}{\h{r}{b}}\int_{|\tau| \delta \leq 1} |\tau|^{k-1}|\widehat{ g}(\tau)|d \tau ,\]
where
\begin{eqnarray*}
\int_{|\tau| \delta \leq 1} |\tau|^{k-1}|\widehat{ g}(\tau)|d \tau & \leq & {\delta}^{1-k}\int_{|\tau| \delta \leq 1}\langle \tau \rangle^{-b'}\langle \tau \rangle^{b'}|\widehat{ g}(\tau)| d \tau \\
& \leq & {\delta}^{1-k} \left(\int_{|\tau| \delta \leq 1}\langle \tau \rangle^{-rb'}d \tau \right)^{\frac{1}{r}}\n{g}{\h{r}{b'}} \\
& \leq & c {\delta}^{\frac{1}{r'}+ b'-k}\n{g}{\h{r}{b'}}
\end{eqnarray*}
and
\begin{eqnarray*}
\|t^k \psi_{\delta}\|^{r'}_{\h{r}{b}} & = & \int \langle \tau \rangle^{r'b} | \widehat{ \psi_{\delta} }^{(k)}(\tau)|^{r'} d \tau \\
& = & {\delta}^{(k+1)r'}\int \langle \tau \rangle^{r'b} |\widehat{ \psi}^{(k)}(\delta \tau)|^{r'} d \tau \\
& \leq & c {\delta}^{(k-b+1)r' - 1}\int \langle \tau \rangle^{r'b}|\widehat{ \psi}^{(k)}(\tau)|^{r'} d \tau =c {\delta}^{(k-b+1)r' - 1} \|t^k \psi\|^{r'}_{\h{r}{b}}.
\end{eqnarray*}
Now a simple computation using the support condition on $\psi$ shows that
\[\|t^k \psi\|_{\h{r}{b}} \le c \|t^k \psi\|_{W^{2,1}} \le c 2^k k^2 \| \psi\|_{W^{2,1}},\]
hence
\[\n{I}{\h{r}{b}} \leq c {\delta}^{1+b'-b} \sum_{k \geq 1} \frac{2^k k^2}{k!} \| \psi\|_{W^{2,1}}\n{g}{\h{r}{b'}}\leq c{\delta}^{1+b'-b}\n{g}{\h{r}{b'}}.\]
Next we consider the second contribution: For $b\geq 0 \ge b' > -\frac{1}{r'}$ we obtain
\begin{eqnarray*}
\n{II}{\h{r}{b}} & \leq & c \n{\psi_{\delta}}{\h{r}{b}}\int_{|\tau| \delta \geq 1} |\tau|^{-1}|\widehat{ g}(\tau)| d \tau \\
& \leq & c {\delta}^{\frac{1}{r} -b}\n{g}{\h{r}{b'}} \left(\int_{|\tau| \delta \geq 1} |\tau|^{-r} \langle \tau \rangle^{-rb'} d\tau \right)^{\frac{1}{r}}\\
& \leq & c {\delta}^{1+b'-b}\n{g}{\h{r}{b'}}.
\end{eqnarray*}
Finally, for the integral $J$ arising in $III$ we have $\widehat{J} = c   \tau^{-1}{\chi}_{|\tau| \delta \geq 1}\widehat{ g }$ and thus
\begin{eqnarray*}
\|J\|^{r'}_{\h{r}{b}}  & \leq & c \int_{|\tau| \delta \geq 1}\langle \tau \rangle^{(b - 1 -b')r'}\langle \tau \rangle^{r'b'}|\widehat{ g}(\tau)|^{r'} d \tau \\
& \leq & c \sup_{|\tau| \geq \frac{1}{\delta}}|\tau|^{(b - 1 -b')r'} \|g\|^{r'}_{\h{r}{b'}}.
\end{eqnarray*}
 For all $b, b' \in {\R}$ satisfying $b-b'\leq 1$ this gives $\n{J}{\h{r}{b}} \leq c {\delta}^{1+b'-b}\n{g}{\h{r}{b'}}$.
For the Fourier transform of the product $\psi_{\delta} J$ we have
\begin{eqnarray*}
\langle \tau \rangle^b |\widehat{ (\psi_{\delta} J)}(\tau)| & \le & \langle \tau \rangle^b \int d \tau_1 |\widehat{ \psi_{\delta}} (\tau_1) \widehat{J}(\tau- \tau_1) | \\
& \le & c \int d \tau_1 |\tau_1|^b |\widehat{ \psi_{\delta}} (\tau_1) \widehat{ J }(\tau- \tau_1)| \\
&+& \int d \tau_1 |\widehat{ \psi_{\delta} }(\tau_1)| \langle \tau- \tau_1 \rangle^b | \widehat{ J} (\tau- \tau_1)|.
\end{eqnarray*}
Hence
\begin{eqnarray*}
\n{\psi_{\delta}J}{\h{r}{b}} & \leq & \n{(|\tau|^b |\widehat{ \psi_{\delta}}|) * |\widehat{ J}|}{L^{r'}_{\tau}} + \n{|\widehat{ \psi_{\delta}}| * (\langle \tau \rangle^b |\widehat{ J}|)}{L^{r'}_{\tau}}\\
& \leq & \n{|\tau|^b |\widehat{ \psi_{\delta}}|}{L^1_{\tau}} \n{\widehat{J}}{L^{r'}_{\tau}} + \n{\widehat{ \psi_{\delta}}}{L^1_{\tau}} \n{J}{\h{r}{b}}\\
& \leq & c (\delta^{-b}\n{J}{\h{r}{0}} +  \n{J}{\h{r}{b}}) \le {\delta}^{1+b'-b}\n{g}{\h{r}{b'}}.
\end{eqnarray*}
Thus (\ref{6}) is shown. Now, if $g$ is a function of both, the time and space variable, it follows that for fixed $\xi$:
\[\int \langle \tau \rangle^{r'b} |\widehat{ Kg}(\xi,\tau)|^{r'} d \tau 
\leq c {\delta}^{r'(1+b'-b)} \int\langle \tau \rangle^{r'b'} |\widehat{ g}(\xi,\tau)|^{r'} d \tau .\]
Multiplying with $\langle \xi \rangle^{r's}$ and integrating with respect to $\xi$ we obtain
\[\|Kg\|^{r'}_{\hh{r}{s}{b}} \leq c  \delta^{r'(1+b'-b)} \|g\|^{r'}_{\hh{r}{s}{b'}}.\]
Applied to $g(t)=U_{\phi} (-t) F(t)$ this gives the desired estimate.
$\hfill \Box$

\subsection{A general local wellposedness theorem}

Here we shall derive a general LWP result for the Cauchy problem
\begin{equation}\label{10}
{\partial}_t u - i \phi (D) u= N(u), \hspace{1cm}u(0)= u_0 \in \h{r}{s},
\end{equation}
where $N$ is a nonlinear function (of degree $\alpha > 1$) of $u$ and its spatial derivatives. Here a solution of (\ref{10}) is understood as a solution of the corresponding integral equation
\begin{equation}\label{11}
u(t)=\Lambda u (t) := U_{\phi} (t) u_0  + {U_{\phi}} *_R N(u) (t).
\end{equation}
For this purpose we introduce the restriction norm spaces
\[\x(\delta) := \{f = \tilde{f}|_{[-\delta,\delta] \times \R^n} : \tilde{f} \in \x\}\]
with norm
\[\n{f}{\x(\delta)}:= \inf \{ \n{\tilde{f}}{\x} : \tilde{f}|_{[-\delta,\delta] \times \R^n} =f\}  .\]
The following theorem reduces the question of local wellposedness completely to nonlinear estimates in $\x$ - spaces:

\begin{satz}[General local wellposedness]\label{t1} Assume that for given $s \in \R$, $r \in (1, \infty)$ there exist $b > \frac{1}{r}$ and $b' \in (b-1,0]$, such that the estimates
\begin{equation}\label{12}
\n{N(u)}{\X{r}{s}{b'}} \le c \|u\|^{\alpha}_{\x}
\end{equation}
and
\begin{equation}\label{13}
\n{N(u) - N(v)}{\X{r}{s}{b'}} \le c (\|u\|^{\alpha -1}_{\x} + \|v\|^{\alpha -1}_{\x})\|u - v\|_{\x}
\end{equation}
are valid. Then there exist $\delta = \delta (\n{u_0}{\h{r}{s}}) > 0$ and a unique solution $u \in \x(\delta)$ of (\ref{10}). This solution is persistent and the mapping data upon solution: $u_0 \mapsto u$, $\h{r}{s} \rightarrow \x (\delta _0)$ is locally Lipschitz continuous for any $\delta _0 \in (0,\delta )$.
\end{satz}

Since the argument is exactly the same as in the $r=2$ - case, it will be sufficient to give a

\vspace{0.3cm}

Sketch of proof: For $u \in \x (\delta)$ with extension $\tilde{u} \in \x$ an extension of $\Lambda u$ is given by
\[\widetilde{\Lambda u}(t)= \psi(t) U_{\phi}(t) u_0 + \psi_{\delta} (t) U_{\phi}*_R N(\tilde{u})(t),\]
where $\psi$ is a smooth cut off function and $\psi_{\delta}(t)=\psi(\frac{t}{\delta})$. Thus
\[\n{\Lambda u}{\x(\delta)} \le \n{\psi U_{\phi} u_0 }{\x} + \n{\psi_{\delta}U_{\phi}*_R N(\tilde{u})}{\x}\]
Using (\ref{5}), Lemma \ref{l2} and (\ref{12}) we see that :
\[\n{\Lambda u}{\x(\delta)} \le c \n{u_0}{\h{r}{s}} + c \delta^{1-b+b'}\|\tilde{u}\|^{\alpha}_{\x}.\]
This holds for all extensions $\tilde{u} \in \x$ of $u \in \x (\delta)$, hence
\[\n{\Lambda u}{\x(\delta)} \le c \n{u_0}{\h{r}{s}} + c \delta^{1-b+b'}\|u\|^{\alpha}_{\x(\delta)}.\]
Similarly, Lemma \ref{l2} and (\ref{13}) give for $u,v \in \x (\delta)$:
\[\n{\Lambda u - \Lambda v}{\x(\delta)} \le c \delta^{1-b+b'}(\|u\|^{\alpha - 1}_{\x(\delta)} + \|v\|^{\alpha - 1}_{\x(\delta)})\|u-v\|_{\x(\delta)} .\]

This shows that for $R=2c\n{u_0}{\h{r}{s}}$ and $\delta^{1-b+b'} \le \frac{1}{4cR^{\alpha - 1}}$ (observe that $1-b+b' > 0$ by assumption) the mapping $\Lambda$ is a contraction of the closed ball of radius $R$ in $\x (\delta)$ into itself. The contraction mapping principle now guarantees the existence of a solution of $\Lambda u = u$. Because of $b>\frac{1}{r}$ any solution $u \in \x (\delta)$ is persistent by (\ref{52}). Thus the standard argument to show uniqueness (see e. g. Prop. 4.2 in \cite{CW90}) in the whole space $ \x (\delta)$ applies. Finally, the statement about continuous dependence can be shown in a straightforward manner using the same estimates as above.
$\hfill \Box$

\vspace{0.3cm}

Here we close the general part of the exposition and specify to the phase function $\phi :\,\,\R \rightarrow \R,\,\,\,\xi \mapsto \xi ^3$, which corresponds to the Airy equation. So, in the sequel the spaces $\x$ are always those defined by this particular phase function.

\section{Bilinear and linear Airy-estimates}

Here we start with the bilinear estimate mentioned in the introduction:
\begin{lemma}\label{l3} 
Let $I^s$ denote the Riesz potential of order $-s$ and let $I^s_-(f,g)$ be defined by its Fourier transform (in the space variable):
\[\widehat{ I_-^s (f,g)} (\xi) := \int_{\xi_1+\xi_2=\xi}d\xi_1|\xi_1-\xi_2|^s \widehat{f}(\xi_1)\widehat{g}(\xi_2).\]
Then we have
\[\n{I^{\frac{1}{2}}I_-^{\frac{1}{2}}(e^{-t\partial ^3}u_1,e^{-t\partial ^3}u_2)}{L^2_{xt}} \le c \n{u_1}{L^2_{x}}\n{u_2}{L^2_{x}}.\]
\end{lemma}

Proof: We will write for short $\int_* d\xi_1$ instead of $\int_{\xi_1+\xi_2=\xi}d\xi_1$. Then, using Fourier-Plancherel in the space variable we obtain:

\begin{eqnarray*}
&&  \q{I^{\frac{1}{2}}I_-^{\frac{1}{2}}(e^{-t\partial ^3}u_1,e^{-t\partial ^3}u_2)}{L^2_{xt}} \\
&=& c\int d\xi |\xi | dt \left| \int_* d\xi_1 |\xi_1-\xi_2|^{\frac{1}{2}} e^{it(\xi_1^3 + \xi_2^3)}\hat{u}_1(\xi_1)\hat{u}_2(\xi_2) \right|^2 \\
&=& c\int d\xi |\xi |dt\int_* d\xi_1 d \eta_1 e^{it(\xi_1^3 + \xi_2^3- \eta_1^3 - \eta_2^3)}(|\xi_1-\xi_2||\eta_1-\eta_2|)^{\frac{1}{2}}\prod_{i=1}^2 \hat{u_i}(\xi_i) \overline{\hat{u_i}(\eta_i)} \\
&=& c\int d\xi |\xi |\int_* d\xi_1 d \eta_1 \delta (\eta_1^3 + \eta_2^3 - \xi_1^3 - \xi_2^3)(|\xi_1-\xi_2||\eta_1-\eta_2|)^{\frac{1}{2}}\prod_{i=1}^2 \hat{u_i}(\xi_i) \overline{\hat{u_i}(\eta_i)} \\
&=& c\int d\xi |\xi |\int_* d\xi_1 d \eta_1 \delta (3\xi (\eta_1^2 - \xi_1^2 + \xi(\xi_1-\eta_1)))(|\xi_1-\xi_2||\eta_1-\eta_2|)^{\frac{1}{2}}\prod_{i=1}^2 \hat{u_i}(\xi_i) \overline{\hat{u_i}(\eta_i)}.
\end{eqnarray*}
Now we use $\delta(g(x)) = \sum_n \frac{1}{|g'(x_n)|} \delta(x-x_n)$, where the sum is taken over all simple zeros of $g$, in our case:
\[g(x)= 3 \xi (x^2 + \xi (\xi_1-x)-\xi_1^2)\]
with the zeros $x_1 =\xi_1$ and $x_2 =\xi-\xi_1$, hence $g'(x_1)=3 \xi(2\xi_1-\xi)$ respectively $g'(x_2)=3 \xi (\xi-2\xi_1)$.
So the last expression is equal to
\begin{eqnarray*}
&&c\int d\xi |\xi | \int_* d\xi_1 d \eta_1\frac{1}{|\xi||2\xi_1-\xi|}\delta(\eta_1-\xi_1)(|\xi_1-\xi_2||\eta_1-\eta_2|)^{\frac{1}{2}}\prod_{i=1}^2 \hat{u_i}(\xi_i) \overline{\hat{u_i}(\eta_i)}\\
&+& c\int d\xi |\xi | \int_* d\xi_1 d \eta_1\frac{1}{|\xi||2\xi_1-\xi|}\delta(\eta_1-(\xi-\xi_1))(|\xi_1-\xi_2||\eta_1-\eta_2|)^{\frac{1}{2}}\prod_{i=1}^2 \hat{u_i}(\xi_i) \overline{\hat{u_i}(\eta_i)}\\
&=&c\int d\xi \int_* d\xi_1 \prod_{i=1}^2|\hat{u_i}(\xi_i)|^2 + c\int d\xi \int_* d\xi_1\hat{u}_1(\xi_1)\overline{\hat{u}_1}(\xi_2)\hat{u}_2(\xi_2)\overline{\hat{u}_2}(\xi_1)\\
&\le& c (\prod_{i=1}^2 \q{u_i}{L^2_x} + \q{\hat{u}_1\hat{u}_2}{L^1_{\xi}}) \le c \prod_{i=1}^2 \q{u_i}{L^2_x}.
\end{eqnarray*}
$\hfill \Box$

\begin{kor}\label{k1} Let $b > \frac{1}{2} \ge s \ge 0$, $\tilde{b}> \frac{1}{6} + \frac{2 s}{3}$. Then the following estimate holds true:
\[\n{I^s I_-^s (u,v)}{L^2_{xt}}\le c \n{u}{\XX{0}{b}}\n{v}{\XX{0}{\tilde{b}}}\]
\end{kor}

Proof: The case $s=\frac{1}{2}$ follows from Lemma \ref{l3}, while in the case $s=0$ we have for $\tilde{b}> \frac{1}{6}$
\[\n{uv}{L^2_{xt}} \le \n{u}{L^8_{xt}} \n{v}{L^{\frac{8}{3}}_{xt}} \le c \n{u}{\XX{0}{b}} \n{v}{\XX{0}{\tilde{b}}}\]
by the well known $L^8_{xt}$ Strichartz type estimate (cf. Corollary \ref{k4} below and its proof) and the trivial case $\XX{0}{0} = L^2_{xt}$. For $0<s<\frac{1}{2}$ we write $w = \Lambda ^{\tilde{b}} v$, where $\Lambda ^{\tilde{b}}$ is defined by $\widehat{\Lambda ^{\tilde{b}} v}(\xi , \tau ) =\langle \tau - \xi ^3 \rangle^{\tilde{b}} \widehat{v} (\xi , \tau )$. Then we have to show that
\begin{equation}\label{100}
\n{I^s I_-^s (u, \Lambda ^{-\tilde{b}}w)}{L^2_{xt}} \le c \n{u}{\XX{0}{b}} \n{w}{L^2_{xt}},
\end{equation}
where
\[\n{I^s I_-^s (u, \Lambda ^{-\tilde{b}}w)}{L^2_{xt}} = \n{|\xi|^s \int_* d \xi _1 d \tau _1|\xi_1 - \xi_2|^s\widehat{u}(\xi_1 , \tau_1 )\langle \tau_2 - \xi_2^3\rangle^{-\tilde{b}}\widehat{w}(\xi_2 , \tau_2 )}{L^2_{\xi \tau}}.\]
By the preceeding (\ref{100}) is already known in the limiting cases $(s,\tilde{b}) = (0, \frac{1}{6}+ \varepsilon)$ and $(s,\tilde{b}) = ( \frac{1}{2}, \frac{1}{2} + \varepsilon)$. Choosing $ \varepsilon = \tilde{b} - \frac{1}{6} - \frac{2 s}{3}$ we have
\[|\xi|^s |\xi_1 - \xi_2|^s \langle \tau_2 - \xi_2^3\rangle^{-\tilde{b}} \le\langle \tau_2 - \xi_2^3\rangle^{-\frac{1}{6}- \varepsilon} + |\xi|^{\frac{1}{2}} |\xi_1 - \xi_2|^{\frac{1}{2}} \langle \tau_2 - \xi_2^3\rangle^{-\frac{1}{2}- \varepsilon}\]
and hence
\[\n{I^s I_-^s (u, \Lambda ^{-\tilde{b}}w)}{L^2_{xt}} \le \n{u \Lambda ^{-\frac{1}{6}- \varepsilon}w}{L^2_{xt}} + \n{I^{\frac{1}{2}} I_-^{\frac{1}{2}} (u, \Lambda ^{-\frac{1}{2}- \varepsilon}w)}{L^2_{xt}} \le c \n{u}{\XX{0}{b}}\n{w}{L^2_{xt}}.\]
$\hfill \Box$

In order to dualize Corollary \ref{k1} we introduce the bilinear operator $I_+^s$ by
\[\widehat{ I_+^s (f,g)} (\xi) := \int_{\xi_1+\xi_2=\xi}d\xi_1|\xi + \xi_2|^s \widehat{f}(\xi_1)\widehat{g}(\xi_2)\]
and the linear operators

\newpage

\[M^s_u v := I_-^s(u,v) \hspace{1cm}\mbox{and}\hspace{1cm}N^s_u w := I_+^s(w,\overline{u}).\]
Then it is easily checked that $M^s_u$ and $N^s_u$ are formally adjoint with respect to the inner product on $L^2_{xt}$. Corollary \ref{k1} expresses the boundedness of
\[M_u^s: \XX{0}{\tilde{b}} \rightarrow L_t^2(\dot{H}^{s})\hspace{1cm}\mbox{with}\hspace{1cm}\|M_u^{s}\| \le c \n{u}{\XX{0}{b}}.\]
But then, of course,
\[N_u^{s}: L_t^2(\dot{H}^{-s}) \rightarrow \XX{0}{-\tilde{b}}\]
is bounded with the same norm, which gives
\begin{kor}\label{k2} Let $b > \frac{1}{2} \ge s \ge 0$, $\tilde{b}> \frac{1}{6} + \frac{2 s}{3}$. Then 
\[\n{I_+^{s} (I^{s}w,u)}{\XX{0}{-\tilde{b}}}\le c \n{w}{L^2_{xt}}\n{u}{\XX{0}{b}}.\]
\end{kor}

Combining Corollary \ref{k2} with the trivial endpoint case of the Hausdorff-Young-inequality we obtain

\begin{kor}\label{k10}
For $1 < r \le 2$, $0 \le \sigma \le \frac{1}{r'}< \beta$ and $b' < -\frac{1}{3}(\frac{1}{r'} + 2 \sigma)$ the following estimate is valid:
\[\n{I_+^{\sigma} (I^{\sigma}w,u)}{\X{r}{0}{b'}}\le c \n{w}{L^2_{xt}}\n{u}{\XX{0}{\beta}}\]
\end{kor}

Proof: Writing $v=\Lambda ^{\beta}u$ we have to show that
\[\n{I_+^{\sigma} (I^{\sigma}w, \Lambda ^{-\beta}v)}{\X{r}{0}{b'}}\le c \n{w}{L^2_{xt}}\n{v}{L^2_{xt}},\]
where
\[\n{I_+^{\sigma} (I^{\sigma}w, \Lambda ^{-\beta}v)}{\X{r}{0}{b'}} = \n{\langle \tau - \xi ^3 \rangle ^{b'} I(\xi , \tau)}{L^{r'}_{\xi \tau}}\]
with
\[I(\xi , \tau) = \int_* d \xi _1 d \tau _1 |\xi + \xi_2|^{\sigma} |\xi_1|^{\sigma}\widehat{w}(\xi_1 , \tau_1 )\langle \tau_2 - \xi_2^3\rangle^{-\beta}\widehat{v}(\xi_2 , \tau_2 ).\]
We choose $\theta = \frac{2}{r}- 1 \in [0,1)$, so that $1-\theta = \frac{2}{r'}$, $s= \frac{\sigma}{1-\theta}$, $\tilde{b} = -\frac{b'}{1-\theta}$ and $b = \frac{\beta}{1-\theta}$. Now H\"older's inequality gives
\begin{eqnarray*}
I(\xi , \tau) & \le & \left(\int_* d \xi _1 d \tau _1 \widehat{w}(\xi_1 , \tau_1 )\widehat{v}(\xi_2 , \tau_2 ) \right)^{\theta} \times \\
&& \left(\int_* d \xi _1 d \tau _1 |\xi + \xi_2|^{s} |\xi_1|^{s}\langle \tau_2 - \xi_2^3\rangle^{-b}\widehat{w}(\xi_1 , \tau_1 )\widehat{v}(\xi_2 , \tau_2 )\right)^{1 - \theta},
\end{eqnarray*}
where the first factor is bounded by
\[\|\widehat{wv}\|^{\theta}_{L^{\infty}_{\xi \tau}} \le \|wv\|^{\theta}_{L^1_{xt}}\le \|w\|^{\theta}_{L^2_{xt}}\|v\|^{\theta}_{L^2_{xt}}.\]
Multiplying the second factor with $\langle \tau - \xi ^3\rangle ^{b'}$ and taking the $L^{r'}_{\xi \tau}$-norm we get the upper bound
\[\|I_+^s (I^s w , \Lambda ^{-b}v)\|^{1-\theta}_{\XX{0}{-\tilde{b}}}.\]
Now the assumptions on $r, \sigma , \beta , b'$ and our choice of $s,b, \tilde{b}$ admit the application of Corollary \ref{k2} and the proof is complete.
$\hfill \Box$

\vspace{0,3cm}

Next we turn to the linear estimates. A slight modification of the proof of Lemma \ref{l3} leads to the following result:

\newpage

\begin{lemma}\label{l4} For $4<q<\infty$ and $\frac{1}{r}=\frac{1}{2}+\frac{1}{q}$ the estimate
\[ \n{e^{-t \partial ^3}u}{L_t^4(\dot{H}^{\frac{1}{4},q})} \le c \n{\widehat{u}}{L^{r'}}\]
is valid.
\end{lemma}

Proof: We assume first that $\widehat{u}= \chi_{[0,\infty)}\widehat{u}$ and write $v = I^{\frac{1}{4}}u$. Then
\begin{eqnarray*}
\|e^{-t \partial ^3}u\|^4_{L_t^4(\dot{H}^{\frac{1}{4},q})} & = & \|e^{-t \partial ^3}v\|^4_{L_t^4(L_x^q)} \\
& = & \| |e^{-t \partial ^3}v|^2 \|^2_{L_t^2(L_x^{\frac{q}{2}})}\\
& \le  & c \| I^{\varepsilon} |e^{-t \partial ^3}v|^2 \|^2_{L^2_{xt}}
\end{eqnarray*}
for $\varepsilon = \frac{1}{2}-\frac{2}{q}$ by Sobolev's embedding theorem. Proceeding as in the proof of Lemma \ref{l3} we get two contributions $I$ and $II$, where - up to constants -
\[ I = \int d\xi |\xi|^{2 \varepsilon -1} \int_* d\xi_1 |\xi_1 -\xi_2|^{-1} |\hat{v}(\xi_1)|^2 |\hat{v}(-\xi_2)|^2 \]
and
\[ II = \int d\xi |\xi|^{2 \varepsilon -1} \int_* d\xi_1 |\xi_1 -\xi_2|^{-1} \hat{v}(\xi_1) \hat{\overline{v}}(\xi_2) \overline{\hat{v}}(\xi_2) \overline{ \hat{\overline{v}}}(\xi_1).\]
By the support condition on $\hat{u}$ resp. $\hat{v}$ we see that the integrand in $I$ is only nonvanishing if $\xi_1 \ge 0 \ge \xi_2$, leading to
\[|\xi_1 -\xi_2|^{-1} =(|\xi_1| +|\xi_2|)^{-1} \le |\xi_1|^{-\frac{1}{2}} |\xi_2|^{-\frac{1}{2}}.\]
This gives
\begin{eqnarray*}
I & \le & \int d\xi |\xi|^{2 \varepsilon -1} \int_* d\xi_1 |\hat{u}(\xi_1)|^2 |\hat{u}(-\xi_2)|^2 \\
& = & \int d\xi_1 |\hat{u}(\xi_1)|^2 \int d\xi |\xi|^{2 \varepsilon -1} |\hat{u}(\xi_1-\xi)|^2 \\
& \le & \n{|\hat{u}|^2}{L^{\frac{r'}{2}}} \n{|\xi|^{2 \varepsilon -1} \ast |\hat{u}|^2}{L^p},
\end{eqnarray*}
where $1= \frac{1}{p} + \frac{2}{r'}$. Now the first factor is nothing but $\q{\hat{u}}{L^{r'}}$, which also controls the second one by the Hardy-Littlewood-Sobolev inequality. Next we observe that the contribution $II$ vanishes by the support assumption on $\hat{u}$ resp. $\hat{v}$. So, in the special case where $\widehat{u}= \chi_{[0,\infty)}\widehat{u}$, the desired estimate is shown. Now, if $\widehat{w}= \chi_{(-\infty,0]}\widehat{w}$ and $u=\overline{w}$, then by $\widehat{\overline{w}}(\xi) = \overline{\widehat{w}}(-\xi)$ we see that $\widehat{u}= \chi_{[0,\infty)}\widehat{u}$. Thus the estimate is valid for $u$. Hence
\begin{eqnarray*}
\|e^{-t \partial ^3}w\|_{L_t^4(\dot{H}^{\frac{1}{4},q})} & = & \|\overline{e^{-t \partial ^3}w}\|_{L_t^4(\dot{H}^{\frac{1}{4},q})} \\
 & = & \|e^{-t \partial ^3}u\|_{L_t^4(\dot{H}^{\frac{1}{4},q})} \\
& \le & c \n{\widehat{u}}{L^{r'}} = c \n{\widehat{w}}{L^{r'}}.
\end{eqnarray*}
Finally the decomposition $u=u_+ + u_-$ with $\widehat{u_+}= \chi_{[0,\infty)}\widehat{u}$ yields the desired result in the general case.
$\hfill \Box$

\vspace{0.3cm}

The endpoint case $(p,q)=(4,\infty)$ is known to be true, too, see Theorem 2.1 in \cite{KPV91}. Next we use interpolation between Lemma \ref{l4}, the conservation of the $L^2$ \nolinebreak -  norm and the trivial estimate
\[\n{e^{-t \partial ^3}u}{L^{\infty}_{xt}} \le c \n{\widehat{u}}{L^{1}}\]
to obtain

\begin{kor}\label{k3} Let $\frac{1}{r} = \frac{2}{p} + \frac{1}{q}$. Then the estimate
\[\n{e^{-t \partial ^3}u}{L_t^p(\dot{H}^{\frac{1}{p},q})} \le c \n{\widehat{u}}{L^{r'}}\]
holds true, if one of the following conditions is fulfilled:
\begin{itemize}
\item[i)] $0 \le \frac{1}{p} \le \frac{1}{4}$, $0 \le \frac{1}{q} < \frac{1}{4}$ or
\item[$$]
\item[ii)] $\frac{1}{4} \le \frac{1}{q} \le \frac{1}{q} + \frac{1}{p} < \frac{1}{2}$ or
\item[$$]
\item[iii)] $(p,q) = (\infty , 2)$.
\end{itemize}
The case $p=q$ is of special interest, here the conditions reduce to $0 \le \frac{1}{p} = \frac{1}{3r} < \frac{1}{4}$ (Airy version of the Fefferman-Stein-estimate).
\end{kor}

\vspace{0.3cm}

The corresponding $\x$ -estimates are gathered in 

\begin{kor}\label{k4} Under the assumptions on $p,q,r$ as in Corollary \ref{k3} and for $b> \frac{1}{r}$ the estimates
\begin{equation}\label{30}
\n{u}{L_t^p(H^{\frac{1}{p},q})} \le c \n{u}{\X{r}{0}{b}}
\end{equation}
and
\begin{equation}\label{31}
\n{u}{\X{r'}{0}{-b}} \le c \n{u}{L_t^{p'}(H^{-\frac{1}{p},q'})}
\end{equation}
are valid, where in (\ref{31}) $\frac{1}{p} + \frac{1}{p'}=\frac{1}{q}+\frac{1}{q'}= \frac{1}{r} + \frac{1}{r'}=1$.
\end{kor}

Proof: Estimate (\ref{30}) with $\dot{H}^{\frac{1}{p},q}$ instead of $H^{\frac{1}{p},q}$ is an immediate consequence of Lemma \ref{l1} and Corollary \ref{k3}. In order to show how to replace the homogeneous space by the inhomogeneous one we may restrict ourselves to the case $p=4$ by interpolation. For that purpose we recall the well known Strichartz type
estimate
\[\n{e^{-t \partial ^3}u}{L^8_{xt}} \le c \n{u}{L_x^2},\]
(which can be obtained from Corollary \ref{k3} by Sobolev's embedding Theorem) respectively its $\xx$ - version
\[\n{u}{L^8_{xt}} \le c \n{u}{\XX{0}{b_0}}, \hspace{1cm}b_0 > \frac{1}{2}.\]
Interpolation with the trivial case $L^2_{xt} = \XX{0}{0}$ gives
\[\n{u}{L^4_{xt}} \le c \n{u}{\XX{0}{b_1}}, \hspace{1cm}b_1 > \frac{1}{3}.\]
Now if $p = \mathcal{F}_x^{-1} \chi_{\{|\xi| \le 1\}}\mathcal{F}_x$ we get for $4<q \le \infty$, $\frac{1}{r}=\frac{1}{2}+\frac{1}{q}$, $b_1 > \frac{1}{3}$ and $b > \frac{1}{r}$:
\[\n{pu}{L_t^4(L_x^q)} \le c \n{J^s pu}{L^4_{xt}} \le c
\n{J^s pu}{\XX{0}{b_1}} \le c \n{u}{\X{r}{0}{b}},\,\,\,(s>\frac{1}{2})\]
where we have used Sobolev's embedding theorem in the space variable as well as the embedding (\ref{4}). This shows (\ref{30}). Finally (\ref{31}) follows from (\ref{30}) by duality.
$\hfill \Box$

\section{The nonlinear estimate}

\begin{satz}\label{t2} Let $2 \ge r > \frac{4}{3}$ and $s \ge s(r)=\frac{1}{2} - \frac{1}{2r}$. Then for all $b' < \frac{1}{2r} -  \frac{5}{8}$ and $b > \frac{1}{r}$ the estimate
\begin{equation}\label{40}
\n{\partial _x (\prod_{i=1}^3 u_i)}{\X{r}{s}{b'}}\le c \prod_{i=1}^3 \n{u_i}{\x}
\end{equation}
holds true.
\end{satz}

Proof: Without loss of generality we may assume that $s = s(r)$. Then we rewrite the left hand side of (\ref{40}) as
\[\n{\langle \tau -  \xi^3 \rangle^{b'} \langle  \xi \rangle^{s}|\xi|\int d \nu  \prod_{i=1}^3 \widehat{u_i}(\xi_i,\tau_i)}{L^{r'}_{\xi,\tau}},\]
where $d \nu = d\xi_1 d\xi_{2} d\tau_1 d \tau_{2}$ and $\sum_{i=1}^3 (\xi_i,\tau_i) = (\xi, \tau)$. We divide the domain of integration into three regions $A$, $B$ and $C$, where in $A$ we assume that $|\xi_1| \sim |\xi_2| \sim |\xi_3|$ or that $|\xi_i| \le 1$ for all $1 \le i \le 3$. Then for all $b' \le 0$ and $b> \frac{1}{r}$ the contribution from this region can be controlled by
\begin{eqnarray*}
\n{\int d \nu  \prod_{i=1}^3 \langle  \xi_i \rangle^{\frac{s+1}{3}} \widehat{u_i}(\xi_i,\tau_i)}{L^{r'}_{\xi,\tau}} & \le & c \n{\prod_{i=1}^3 J^{\frac{s+1}{3}}u_i}{L^r_{xt}}\\
\le c \prod_{i=1}^3 \n{J^{\frac{s+1}{3}}u_i}{L^{3r}_{xt}} & \le & c \prod_{i=1}^3 \n{u_i}{\x} ,
\end{eqnarray*}
where we have used the Hausdorff-Young and H\"older inequalities as well as Corollary \ref{k4}. In region $B$ we assume that $|\xi _{max}| \sim |\xi _{med}| \gg |\xi _{min}|$, where $\xi _{max}$, $\xi _{med}$ and $\xi _{min}$ are defined by $|\xi _{max}| \ge |\xi _{med}| \ge |\xi _{min}|$. By symmetry it is sufficient to consider the subregion, where $|\xi_1| \ge |\xi_2| \ge |\xi_3|$. The contribution of this subregion is bounded by
\begin{equation}\label{41}
c \n{(J^{s + \frac{1}{4} - \mu}u_1)I^{\frac{1}{2}}I_-^{\frac{1}{2}}(J^{-\frac{1}{8}}u_2,J^{-\frac{1}{8}}u_3)}{\X{r}{\mu}{b'}},
\end{equation}
where $\mu = \frac{1}{4} - \frac{1}{3r}$. By interpolation between the $p=q$ - part of (\ref{31}) and Hausdorff-Young (with constant $r$) we obtain
\[\n{u}{\X{r}{\mu}{b'}} \le c \n{u}{L^q_{xt}},\]
whenever
\[\mu = \frac{1}{2}(\frac{1}{q}-\frac{1}{r}) \,\,\, , \hspace{1cm} b' < -3\mu \,\,\, , \hspace{1cm}\frac{1}{2} \le \frac{1}{r}\le \frac{1}{q} \le \frac{2}{3}+ \frac{1}{3r}.\]
Applying this to (\ref{41}) with $\frac{1}{q} = \frac{1}{2} + \frac{1}{3r}$ (i. e. $\mu = \frac{1}{4} - \frac{1}{3r}$ as chosen above) we get - for all $b' < \frac{1}{r} - \frac{3}{4}$ - the upper bound
\begin{eqnarray*}
&& c \n{(J^{s + \frac{1}{4} - \mu}u_1)I^{\frac{1}{2}}I_-^{\frac{1}{2}}(J^{-\frac{1}{8}}u_2,J^{-\frac{1}{8}}u_3)}{L^q_{xt}} \\
&& \le c \n{J^{s + \frac{1}{4} - \mu}u_1}{L^{3r}_{xt}} \n{I^{\frac{1}{2}}I_-^{\frac{1}{2}}(J^{-\frac{1}{8}}u_2,J^{-\frac{1}{8}}u_3)}{L^2_{xt}}\\
&& \le c \n{J^{ \frac{1}{4} - \mu - \frac{1}{3r}}u_1}{\x} \n{u_2}{\XX{-\frac{1}{8}}{\frac{1}{2}+}}\n{u_3}{\XX{-\frac{1}{8}}{\frac{1}{2}+}}
\end{eqnarray*}
by H\"older's inequality, (\ref{30}) and Corollary \ref{k1}. The first factor is nothing but $\n{u_1}{\x}$ by our choice of $\mu$. For the second and third factor we use the embedding (\ref{4}) to obtain
\[\n{u}{\XX{-\frac{1}{8}}{\frac{1}{2}+}} \le c \n{u}{\x},\]
whenever $b > \frac{1}{r}$ and $s + \frac{1}{8} > \frac{1}{r} - \frac{1}{2}$. The latter condition is fulfilled since $r>\frac{4}{3}$. Next we consider the region $C$, where $|\xi_{max}| \gg |\xi_{med}|$. Again we may restrict ourselves to the subregion with $|\xi_1| \ge |\xi_2| \ge |\xi_3|$. Using the notation $f \preceq g$ for $|\widehat{f}| \le c |\widehat{g}|$, we have that in this subregion
\begin{eqnarray*}
&& J^s \partial_x (u_1u_2u_3) \\
& \preceq & J^s (I^{\frac{1}{2}}I_-^{\frac{1}{2}}(u_1,u_2)u_3) \hspace{2cm} (\mbox{since}\,\,\,|\xi_1 +\xi_2 + \xi_3| \le c |\xi_1 + \xi_2|^{\frac{1}{2}}|\xi_1 - \xi_2|^{\frac{1}{2}})\\
& \preceq & J^{s+\frac{3}{8}} (I^{\frac{1}{2}}I_-^{\frac{1}{2}}(J^{-\frac{1}{8}}u_1,J^{-\frac{1}{8}}u_2)J^{-\frac{1}{8}}u_3) \hspace{1.9cm}(\mbox{since}\,\,\,|\xi_i| \le c |\xi|,\,\,\,1\le i \le 3) \\
& \preceq & I_+^{\sigma}(I^{\sigma}I^{\frac{1}{2}}I_-^{\frac{1}{2}}(J^{-\frac{1}{8}}u_1,J^{-\frac{1}{8}}u_2),J^{-\frac{1}{8}}u_3),\hspace{3.9cm}(\sigma =  \frac{s}{2} + \frac{3}{16})
\end{eqnarray*}
since $|\xi_1 +\xi_2 + \xi_3| \le c |\xi_1 + \xi_2|^{\frac{1}{2}}|\xi_1 + \xi_2 + 2\xi_3|^{\frac{1}{2}}$. Thus the contribution of this region is bounded by
\begin{eqnarray*}
&& \n{I_+^{\sigma}(I^{\sigma}I^{\frac{1}{2}}I_-^{\frac{1}{2}}(J^{-\frac{1}{8}}u_1,J^{-\frac{1}{8}}u_2),J^{-\frac{1}{8}}u_3)}{\X{r}{0}{b'}} \\
& \le & \n{I^{\frac{1}{2}}I_-^{\frac{1}{2}}(J^{-\frac{1}{8}}u_1,J^{-\frac{1}{8}}u_2)}{L^2_{xt}}\n{J^{-\frac{1}{8}}u_3}{\XX{0}{\beta}}
\end{eqnarray*}
by Corollary \ref{k10}, provided $0 \le \frac{s}{2} + \frac{3}{16} \le \frac{1}{r'}<\beta $ (which is fulfilled for any $\beta > \frac{1}{2}$ and for $s=s(r)$, since $2 \ge r > \frac{4}{3}$) and $b' < -\frac{1}{3}(\frac{1}{r'} + s + \frac{3}{8})=\frac{1}{2r} - \frac{5}{8}$ as demanded. Applying Corollary \ref{k1} to the first factor we see that the whole expression is bounded by
\[c \prod_{i=1}^3\n{u_i}{\XX{-\frac{1}{8}}{\frac{1}{2}+}} \le c \prod_{i=1}^3\n{u_i}{\x},\]
where in the last step we have used the embedding (\ref{4}) again. Finally, comparing the restrictions on $b'$ arising from the three regions, we see that in the allowed range for $r$ the last one, i.e. $b' < \frac{1}{2r} - \frac{5}{8}$, is the strongest.
$\hfill \Box$

\vspace{0.3cm}

Choosing $b' \in (- \frac{1}{r'},\frac{1}{2r}-\frac{5}{8})$ and $b \in (\frac{1}{r}, b'+1)$ we see that the assumptions of Theorem \ref{t1} are fulfilled. So the final result is

\begin{satz}\label{t3} For $2 \ge r > \frac{4}{3}$, $s \ge s(r)=\frac{1}{2} - \frac{1}{2r}$ and data $u_0 \in \h{r}{s}$ the Cauchy problem (\ref{1}) is locally well posed in the sense of Theorem \ref{t1}.
\end{satz}

\end{document}